\documentclass[12pt]{article}
\usepackage[utf8]{inputenc}
\usepackage{amsmath}
\usepackage{amssymb}
\usepackage{amsthm}
\newcommand{\h}{\hspace{1cm}}
\newcommand{\z}{\varepsilon}
\newcounter{local}
\newcounter{locallocal}
\newcommand{\scl}{\stepcounter{local}}
\begin{document}
\begin{center}

{\bf On a Reaction-Diffusion System Modeling Infectious Diseases Without Life-time Immunity}

\bigskip
Hong-Ming Yin\footnote{Corresponding Author. Email: hyin@wsu.edu}\\
Department of Mathematics and Statistics\\
Washington State University,\\
Pullman, WA 99164, USA.
\end{center}

\begin{abstract}
    In this paper we study a mathematical model for an infectious disease such as 
    Cholera without life-time immunity. Due to the different mobility for susceptible, infected human and recovered human hosts, the diffusion coefficients are assumed to be different. The resulting system is governed by a strongly coupled reaction-diffusion system with different diffusion coefficients. Global existence and uniqueness are established under certain assumptions on known data. Moreover, global asymptotic behavior of the solution is obtained when some parameters satisfy certain conditions. These results extend the existing results in the literature. The main tool used in this paper comes from the delicate theory of  elliptic and parabolic equations. Moreover, the energy method and Sobolev embedding are used in deriving {\em apriori} estimates. The analysis developed in this paper can be employed to study other epidemic models in biological and ecological systems.

\end{abstract}

\ \\
{\bf Key Words and Phases:} Infectious disease model for Cholera; Reaction-diffusion system, Global existence and Uniqueness; Asymptotic behavior and global attractor.

\newpage
\section{Introduction}

In an ecological environment various infectious diseases in human or animals occur frequently (see \cite{WHO}). The current global Covid-19 pandemic is one such an example. Other recent examples include epidemics caused by the HIV, Cholera and Zika viruses. Scientists often use the well-known Susceptibility-Infection-Recovery (SIR) model (\cite{KM1927}) to describe how a virus spreads and evolves in future. The SIR model and its various extensions have been studied by many scientists ( see \cite{AM1979,MA1979,DW2002,LGWK1999,KP1994,SD2012,TH1992}, for examples). A good review for the model written by Hethcote in 2000 can be seen in \cite{H2000}. The most of these studies focus on  understanding the complicated dynamics of interaction between different hosts and viruses. However, this SIR model and its extensions cannot include the mobility of human hosts around  different geographical regions. In order to take the movement of human hosts into consideration,  one must develop a new mathematical model which can reflect these factors (see \cite{ESTD2002,LN1996,YCW2017}). Toward this goal, considerable progress has been made for different types of infectious diseases. In particular, many researchers have studied the mathematical model for the epidemic caused by the Cholera virus (see \cite{AB2011,DW2002,LGWK1999,LW2011,TW2011,TH1992,WW2015,YW2016,YW2017}, etc.).
 
 In this paper we consider the mathematical model for the cholera epidemic with diffusion processes. A novel feature for the model is that there is no life-time immunity. This leads to a coupled reaction-diffusion system with different diffusion coefficients. We begin by describing the model system recently studied in \cite{TW2011,WW2015,YW2017,Y2018}.
 
 Let $\Omega$ be a bounded domain in $R^n$ with $C^2$-boundary $\partial \Omega$. 
 Let $Q_T=\Omega \times (0,T]$ for any $T>0$. When $T=\infty$, we denote $Q_{\infty}$ by $Q$.
 Let $S(x,t), I(x,t)$ and $ R(x,t)$ represent, respectively, susceptible, infected and recovered human hosts. Let $B(x,t)$ be the concentration of bacteria. 
Then by the population growth and the conservation laws we see that $S,I,R$ and $B$ satisfy the following reaction-diffusion system:
\setcounter{section}{1}
\setcounter{local}{1}
 \begin{eqnarray}
 S_t-\nabla[d_1\nabla S] & = & b(x,t)-\beta_1SI-\beta_2S \cdot h(B)-dS+\sigma R,\\
 I_t-\nabla[d_2\nabla I] & = & \beta_1SI+\beta_2S \cdot h(B)-(d+\gamma)I,\scl \\
 R_t-\nabla[d_3\nabla R] & = & \gamma I-(d+\sigma)R,\scl \\
 B_t-\nabla[d_4\nabla B] & = & -\sum_{k=1}^{n}b_kB_{x_{k}}+\xi I +gB(1-\frac{B}{K})-\delta B\scl
 \end{eqnarray}
 subject to the following initial and boundary conditions:
 \begin{eqnarray}
 & & (\nabla_{\nu}S,\nabla_{\nu}I,\nabla_{\nu}R, \nabla_{\nu}B) = 0, \h (x,t)\in \partial \Omega\times (0,T],\scl\\
 & & (S(x,0),I(x,0),R(x,0),B(x,0))=(S_0(x),I_0(x),R_0(x),B_0(x)),  x\in \Omega,\scl
 \end{eqnarray}
where $\nu$ represents the outward unit normal on $\partial \Omega$, $h(B)$ is a differential function  with \[ h(0)=0, 0\leq h(B)\leq 1.\]
 
A typical example used in \cite{YW2017} for $h(B)$ is
\[ h(B)=\frac{B}{B+K}.\]
For reference, we list various parameters in the model as in \cite{WW2015,YW2016}.
 \begin{eqnarray*}
 d & = & \mbox{ the natural death rate}\\
 \gamma & = & \mbox{the recovery of infectious individuals}\\
 b & = & \mbox{the influx of susceptible host}\\
 \sigma & = & \mbox{the rate at which recovered individuals lose immunity}\\
 \delta & = & \mbox{the natural death rate of bacteria}\\
 \xi & = & \mbox{the shedding rate of bacteria by infectious human hosts},\\
 d_i & = & \mbox{the diffusion coefficients, $i=1,2,3,4$},\\
 b_{k} & = & \mbox{the coefficient of the bacteria convection},\\
 K & = & \mbox{ the maximum capacity of the bacteria}
 \end{eqnarray*}
 
 We would like to point out that the analysis for the reaction-diffusion system (1.1)-(1.4) is often difficult since there is no comparison principle (\cite{Lieberman,PAO}). Some basic questions such as global existence and uniqueness for the system are extremely challenging since Eq.(2.1) and Eq.(2.2) contain some quadratic terms in the system (\cite{BH,CC2017,CDV2009}). Many of the quantitative properties are still open. For example, from the physical point of view, the concentration, $S(x,t),I(x,t), R(x,t)$ and $B(x,t)$  must be nonnegative and bounded in $Q_T$. However, it appears that there has been no rigorous proof in the previous research when the space dimension is greater than 1. 
 Answering these open questions is one of the motivations for this paper.  Moreover, with different diffusion coefficients in a reaction-diffusion system, the dynamics of a solution may be very different from that of an ODE system. The most striking example is the Turing phenomenon in which the solution of an ODE system is stable while the solution of the corresponding reaction-diffusion system is unstable when one diffusion coefficient is much larger than the other (\cite{EVANS}).
 
 When the space dimension is equal to $1$, the authors of \cite{LW2011,SD2012,WW2015,YW2017,Y2018} studied the system (1.1)-(1.6). The global well-posedness is established for the model. Moreover, some global dynamical analysis for the solution is carried out in these papers. However, when the spatial dimension $n$ is greater than 1, the analysis becomes much more complicated. The purpose of this paper is to study the reaction-diffusion system (1.1)-(1.6) when the space dimension is greater than 1.  Moreover, the diffusion coefficients in the model for susceptible, infected and recovered human hosts are different. Furthermore, the diffusion coefficients in our model allow to depend on space and time, which are certainly more realistic than in the previous model. This fact is observed from the recent Covid-19 pandemic in which the mobility for infected patients is close to zero due to required global quarantine and travel restriction. It is also clear that for flu-like viruses the rate of infection is much higher in the winter than in the summer. 
 
 In this paper we establish a global existence result for any space dimension $n$. Particularly, we prove that the solution is actually classical if the dimension is less than or equal to $3$ with no restriction on all parameters. Hence, global well-posedness for the model problem (1.1)-(1.6) is established when the space dimension is less than or equal to $3$. Moreover, under a condition on some parameters, we are able to prove that the solution is uniformly bounded and converges to the steady-state solution (global attractor, see \cite{Temam}) for any space dimension as time evolves. These results improve the previous research obtained by others where the well-posedness is proved when the space dimension is equal to 1 (see \cite{SD2012,YW2016,YW2017,Y2018}).  The main idea for establishing global existence is to derive various {\em apriori} estimates (see \cite{YCW2017}). Our analysis in this paper uses a lot of very delicate results for elliptic and parabolic equations (\cite{BG1989,Lieberman,LSU,Y1997}). Particularly, we use a subtle form of Galiardo-Nirenberg's inequality to derive a uniform bound in $Q$ for the solution when the space dimension is less than or equal to $2$ without any restriction on parameters. The global asymptotic analysis is based on accurate energy estimates for the solution of the system. The method developed in this paper can also be used to deal with different epidemic models caused by some viruses such as avian influenza for birds (\cite{VWZ2012}).

  The paper is organized as follows. In Section 2 we first recall some basic  function spaces
   and then state the main results. In Section 3, we prove the first part of the main results on global solvability of the system (1.1)-(1.6) (Theorem 2.1 and Theorem 2.2). The long-time behavior of the solution in the second part is proved in Section 4 (Theorem 2.3 and Theorem 2.4). In Section 5, we give some concluding remarks.

 \section{ Notation and the Statement of Main Results}
 
 For reader's convenience, we recall some basic Sobolev spaces which are standard in dealing with elliptic and parabolic partial differential equations.
 
 Let $\alpha\in (0,1)$. We denote by $C^{\alpha, \frac{\alpha}{2}}(\bar{Q}_{T})$ the  H\"older space in which every function is H\"older continuous with respect to $(x,t)$ with exponent $(\alpha,\frac{\alpha}{2})$ in $\bar{Q}_{T}$.
 
 Let $p>1$ and $V$ be a Banach space with norm $||\cdot||_v$, we denote
 \[ L^p(0,T; V)=\{ F(t): t\in [0,T]\rightarrow V: ||F||_{L^p(0,T;V)}<\infty\},\]
 equipped with norm
 \[ ||F||_{L^{p}(0,T;V)}=\left(\int_{0}^{T} ||F||_v dt\right)^{\frac{1}{p}}.\]
 When $V=L^{p}(\Omega)$, we simply use
 \[ L^p(Q_{T})=L^{p}(0,T;L^{p}(\Omega)).\]
 Moreover, the $L^p(Q_{T})$-norm is denoted by $||\cdot||_p$ for simplicity.
 
Sobolev spaces $W_{p}^{k}(\Omega)$ and $W_{p}^{k,l}(Q_{T})$ are defined the same as in the classical books such as \cite{EVANS} and \cite{LSU}.

Let $V_2(Q_T)=\{ u(x,t)\in C([0,T];W_{2}^{1,0}(\Omega)): ||u||_{V_{2}}<\infty\} $
 equipped with the norm
 \[ ||u||_{V_{2}}=\sup_{0\leq t\leq T}||u||_{L^{2}(\Omega)}+\sum_{i=1}^{n}||u_{x_{i}}||_{2}.\]
 
 We first state the basic assumptions for the diffusion coefficients and known data.
 All other parameters in Eq.(1.1)-(1.4) are assumed to be positive automatically throughout this paper.
 Since there is no essential difference for the analysis in this paper, 
 we set $b_k=0$ in Eq.(1.4) for all $k$ and choose $h(s)=\frac{s}{s+K}$.
 
 \ \\
 {\bf H(2.1).}  Assume that $d_i(x,t)\in L^{\infty}(\Omega\times (0,\infty))$ for all $i$. There exist two positive constants $d_0$ and $D_0$ such that 
 \[ 0<d_0\leq d_i(x,t)\leq D_0, \h (x,t)\in Q_{T}.\]
 {\bf H(2.2).} Assume that all initial data $S_0(x), I_0(x), R_0(x), B_0(x)$ are nonnegative on $\Omega$. Moreover, $Z_0(x)=(S_0(x),I_0(x),R_0(x),B_0(x))\in C^{\alpha}(\bar{\Omega})^4.$
 
 \  \\
 {\bf H(2.3)}. Let $0\leq b(x,t)\in L^{\infty}(Q)$. Moreover,
 \[ ||b||_{L^{\infty}(Q)}\leq b_0<\infty.\]
 
 For brevity, we set 
 \[ u_1(x,t)=S(x,t), u_2(x,t)=I(x,t), u_3=R(x,t), u_4(x,t)=B(x,t), \, (x,t)\in Q. \]
 We use $Z(x,t)=(u_1,u_2,u_3,u_4)$ to be a vector function defined in $Q$. The right-hand sides of the equations (1.1), (1.2), (1.3) and (1.4) are denoted by $f_1(Z), f_2(Z), f_3(Z),f_4(Z)$, respectively.
 With the new notation, the system (1.1)-(1.6) can be written as the following reaction-diffusion system:
 \setcounter{section}{2}
 \setcounter{local}{1}
 \begin{eqnarray}
 & & u_{1t}-\nabla[d_1\nabla u_1]=f_1(Z), \h (x,t)\in Q_{T},\\
 & & u_{2t}-\nabla[d_2\nabla u_2]=f_2(Z), \h (x,t)\in Q_{T},\scl \\
& &  u_{3t}-\nabla[d_3\nabla u_3]=f_3(Z), \h (x,t)\in Q_{T},\scl\\
 & & u_{4t}-\nabla[d_4\nabla u_4]=f_4(Z), \h (x,t)\in Q_T,\scl 
 \end{eqnarray}
 subject to the initial and boundary conditions:
 \begin{eqnarray}
 & & Z(x,0)=(S_0(x),I_0(x),R_0(x),B_0(x)), \h x\in \Omega, \scl \\
 & & \nabla_{\nu}Z(x,t)=0,           \h (x,t)\in \partial \Omega\times (0,T].\scl
 \end{eqnarray}
 
 Let $q>1$. We define a product space $X$ equipped with the standard product norm:
 \[ X=L^{q}(0,T;W^{1, q}(\Omega))^{4}.\]
 The conjugate dual space of $X$ is denoted by $X^*$.
 \ \\
 {\bf Definition 2.1} We say $Z(x,t)$ to be a weak solution to the problem (2.1)-(2.6) in $Q_{T}$ if 
 $Z(x,t)\in X$ such that
 \begin{eqnarray}
 & & \int_{0}^{T}\int_{\Omega}[-u_k\cdot \phi_{kt}+d_k\nabla u_k\cdot \nabla\phi_k]dxdt=\int_{\Omega} f_k\phi _k(x,0)dx, \scl
 \end{eqnarray}
 for any test function $\phi_k, \phi_{k t}\in X^*$ with $\phi_k(x,T)=0$ on $\Omega$ for all $k=1,2,3,4$.

 \ \\
 {\bf Theorem 2.1.} Under the assumptions H(2.1)-H(2.3), the problem (2.1)-(2.6) possesses at least one  weak solution and the solution is nonnegative in $Q_T$. Moreover, $u_3(x,t)\in W_{p}^{2,1}(Q_{T})$
 for any $p\in (1, \frac{n+2}{n+1})$. 
 
 A much better result can be proved when the space dimension is less than or equal to $2$.
 
 \ \\
 {\bf Corollary 2.1.} Under the assumption H(2.1)-H(2.3)  the weak solution is unique, uniformly bounded and H\"older continuous in $\bar{Q}$ if the space dimension is less than or equal to $2$.
 
 With more regularity for the coefficients and other known functions, we can obtain the classical solution.
 
 \ \\
 {\bf H(2.4).}  Assume that $d_i(x,t)\in C^{1+\alpha, \frac{1+\alpha}{2}}(\bar{Q})$ for $i=1, 2, 3, 4$, and $b(x,t)\in C^{\alpha, \frac{\alpha}{2}}(\bar{Q})$.
  Let $Z_0(x)\in C^{2+\alpha}(\bar{\Omega})$.
 \ \\
 {\bf Theorem 2.2} Under the assumptions H(2.1)-H(2.4), the problem (2.1)-(2.6) has a unique classical solution $Z(x,t)\in C(\bar{Q}_{T})\bigcap C^{2+\alpha, 1+\frac{\alpha}{2}}(Q_{T})$ for any $T>0$ if the space dimension $n\leq 3$. 
 
 With a certain condition for some parameters, we can prove that the solution exists globally and is uniformly bounded in $Q$ for any dimension $n$.
 
 \ \\
 {\bf Theorem 2.3} Under the assumption H(2.1)-(2.4), the problem (2.1)-(2.6) has a unique global solution  in $V_2(Q)^4$ if 
 \[ d-g_0>0,\]
 where \[ g_0=\frac{\sigma+\beta_1+\beta_2+\gamma}{4}.\]
 Moreover, $Y(t)$ is uniformly bounded in $\Omega$, where
 \[ Y(t)=\int_{\Omega}[S^2+I^2+R^2+B^2]dx.\]
 Moreover, the weak solution is smooth, provide all coefficients $d_k$ and $b(x,t),Z_0(x)$ are smooth.
 
 With the result of Theorem 2.3, we can establish the global asymptotic behavior of the solution.
 
 \ \\
 {\bf H(2.5).}  Assume 
 \[\lim_{t\rightarrow \infty}b(x,t)=b_0(x), \lim_{t\rightarrow \infty}d_i(x,t)=d_{i\infty}(x), \h x\in \Omega, i=1, 2, 3, 4,\]
 in $L^2(Q)$-sense. Moreover, $d_{i\infty}(x)\in C^{\alpha}(\bar{\Omega})$.

 \ \\
 {\bf Theorem 2.4.} Under the assumption H(2.1)-H(2.5), if $ d-g_0>0$ and $n\leq 4,$
 then,
 \[ \lim_{t\rightarrow \infty}(S(x,t), I(x,t),R(x,t), B(x,t))=(S_{\infty}(x), 0, 0, 0),\]
 in $L^2(Q)$-sense, where $S_{\infty}(x)$ is the steady-state solution for the following elliptic equation
 \begin{eqnarray*}
 & & -\nabla[d_1(x)\nabla S_{\infty}]=b_0(x)-dS_{\infty}, \h x\in \Omega\\
 & &  \nabla_{\nu}S_{\infty}(x)=0, \h x\in \partial \Omega.
 \end{eqnarray*}
 That is $(S_{\infty}(x), 0, 0, 0)$ is a global attractor.
 
 \ \\
 {\bf Remark 2.1.} If the space dimension $n$ is less than or equal $2$, the condition 
 \[ d-g_0>0\] can be removed by Theorem 2.2.  Whether or not this condition in Theorem 2.4 is necessary is an open question when the space dimension $n\geq 3$.
 
 \section{Global Solvability and Regularity}
 
From the physical model of view, the concentration must be nonnegative. However, this fact has not been proved rigorously in the literature. Here we first prove this fact by using the strong maximum principle. The proof is rather technical. 
\setcounter{section}{3}
\setcounter{local}{1}
 \ \\
 {\bf Lemma 3.1}:  Under the assumptions H(2.1)-H(2.2), each concentration is nonnegative in $Q$,
 \[ u_1(x,t), \geq 0, u_2(x,t)\geq 0, u_3(x,t)\geq 0, u_4(x,t)\geq 0, \h (x,t)\in Q_T.\]
{\bf Proof}. First of all, we assume that all initial functions are positive over $\bar{\Omega}$. Namely, we assume that there exists 
a small number $a_0>0$ such that
\[ \inf_{x\in \Omega}\{ S_0(x), I_0(x), R_0(x), B_0(x)\} \geq a_0>0.\]
By the continuity, we know that there is a number $T_0>0$ such that
\[ u_1(x,t), u_2(x,t), u_3(x,t), u_4(x,t) \geq \frac{a_{0}}{2}, \h (x,t)\in Q_{T_{0}}.\]
Let 
\[ T^*=\sup_{0<t<T_{0}} \{T_0\in [0,\infty): 
u_1(x,t), u_2(x,t), u_3(x,t), u_4(x,t) >0, x\in \Omega.\}\]
If $T^*=\infty$, then the conclusion holds. If $T^*<\infty$, then at least one of quantities attains $0$ at some point $(x^*,T^*)$ for some $x^*\in \bar{\Omega}$. Suppose 
\[ u_1(x^*,T^*)=0.\]
Since $b(x,t)\geq 0,  u_3(x,t)\geq 0$ on $Q_{T^{*}}$, from the boundary condition we can use Hofp's lemma to conclude that
$x^*$ must not be located on the boundary $\partial \Omega$. On the other hand, the strong maximum principle implies that
$x^*$ can not be located at an interior point of $\Omega$. It follows that
\[ \min_{ (x,t)\in Q_{T^*}}u_1(x,t)>0,\]
which is a contradiction.

From Eq. (2.2), we note that
\[ \beta_1u_1u_2\geq 0, \beta_2 u_1 \frac{u_{4}}{u_{4}+K}\geq 0, (x,t)\in Q_{T^*}.\]
We can use the strong maximum principle to conclude that 
\[\min_{ (x,t)\in Q_{T^*}}I(x,t) >0.\]
Next since $\xi>0, \gamma>0$ and $u_2(x,t)\geq 0$ in $ Q_{T^{*}}$, the same argument holds for $u_3(x,t)$
by Eq.(2.3) and $u_4(x,t)$ by Eq.(2.4). Namely,
\[ \min_{ (x,t)\in Q_{T^*}}u_3(x,t)>0;\, \min_{ (x,t)\in Q_{T^*}}u_4(x,t)>0,\]
which is a contradiction with the assumption. This implies that $T^*=\infty$.

Now when initial functions do not have a positive lower bound, for any $\varepsilon>0$ we simply use 
$(S_0(x)+\varepsilon, I_0(x)+\varepsilon, R_0(x)+\varepsilon, B_0(x)+\varepsilon)$ to replace the original initial vector function $Z(x,0)$.
Then we know that the corresponding solution $(u_{1\varepsilon}(x,t),u_{2\varepsilon}(x,t),u_{3\varepsilon}(x,t)$ and $u_{4\varepsilon}(x,t)$
are positive in $Q_T$ for any $T>0$. By taking the limit as $\varepsilon\rightarrow 0$, we obtain the nonnegetivity of $u_k(x,t)$ in $Q$ for $k=1,2,3,4$.

As in the standard analysis in deriving an apriori estimate for solution of a partial differential equation we may assume that the solution is smooth in $Q_{T}$, A special attention is paid to be what a constant $C$ depends precisely on known data in the derivation. We will denote by $C, C_1, C_2\cdots,$ etc.,the generic constants in the derivation. Those constants may be different from one line to  the next as long as their dependence is the same.

\ \\
{\bf Lemma 3.2.} Under the assumptions H(2.1)-H(2.2), there exists a constant $C_1$ such that
\begin{eqnarray*}
& & \sup_{0<t<\infty}\sum_{k=1}^{4}||u_k||_{L^1(\Omega)}\leq C_1,\\
& & \int\int_{Q_{T}}[u_1u_2 +u_1h(u_4)] dxdt\leq C_2,
\end{eqnarray*}
where $C_1$ and $C_2$ depends only on know data. In addition $C_2$ also depends on the upper bound of $T$.\\
{\bf Proof}: We take integration over $\Omega$ for Eq.(2.1) to Eq.(2.3) and then add up to see
\begin{eqnarray*}
& & \frac{d}{dt}\int_{\Omega}(u_1+u_2+u_3)dx+d\int_{\Omega}(u_1+u_2+u_3)dx \\
& & =\int_{\Omega}b(x,t)dx \leq b_0|\Omega|.
\end{eqnarray*}

Since $d>0$, it follows that
\[ \sup_{0<t<\infty}\int_{\Omega}(u_1+u_2+u_3)dx\leq C_1,\]
where $C_1$ depends only on known data.

Now we integrate of $\Omega$ for Eq.(2.1) again to find
\[ \frac{d}{dt}\int_{\Omega}u_1dx+\beta_1\int_{\Omega}u_1u_2 dx+\beta_2\int_{\Omega}u_1h(u_4)dx
= \int_{\Omega}b(x,t)dx +\sigma\int_{\Omega}u_3dx.\]
It follows that
\[ \beta_1\int_{0}^{T}\int_{\Omega}u_1u_2 dxdt+\beta_2\int_{0}^{T}\int_{\Omega}u_1(u_4)dxdt\leq C_2.\]

Next we take integration over $\Omega$ for Eq.(2.4) to obtain
\begin{eqnarray*}
& & \frac{d}{dt}\int_{\Omega}u_4dx+\delta\int_{\Omega}u_4dx+g\int_{\Omega}u_{4}^2dx\\
& & =\xi \int_{\Omega}u_2dx+g\int_{\Omega}u_4 dx\\
& & \leq \xi \int_{\Omega}u_2dx+ \frac{g}{2}\int_{\Omega}u_{4}^2dx+C.
\end{eqnarray*}
Since $L^1(\Omega)$-norm of $u_2$ is uniformly bounded, it follows that
\[\frac{d}{dt}\int_{\Omega}u_4dx+\delta\int_{\Omega}u_4dx\leq C.\]
Thus, we obtain the desired $L^1(\Omega)$-estimate for $u_4$.

\hfill Q.E.D.

To derive more apriori estimates, we need to use a very delicate result from the theory of parabolic equations with measure data. The reader can find the proof from \cite{LSU} for $n\leq 2$ and \cite{BG1989} for $n\geq 3$.

\ \\
{\bf Lemma 3.3.} Let $a(x,t)$ be a measurable function with $0<a_0\leq a(x,t)\leq a_1<\infty$ in $Q$ and $f(x,t)\in L^{\infty}(0,\infty; L^{1}(\Omega))$.
Let $u(x,t)$ be a weak solution of the following parabolic equation
\begin{eqnarray*}
& & u_t-\nabla [a(x,t) \nabla u]=f(x,t), \h (x,t)\in Q,\\
& & \nabla_{\nu} u(x,t)=0, \h (x,t)\in \partial \Omega \times (0,\infty),\\
& & u(x,0)=u_0(x), \h x\in \Omega.
\end{eqnarray*}
Then, 
(a) If $n=1$, then $u(x,t)$ is uniformly bounded in $Q$. Moreover,
\[ \sup_{0<t<\infty}||u||_{\infty}\leq C_3;\]
(b) If $n=2$, then $u(x,t)\in L^p(Q)$ for any $p\in (1, \infty)$. Moreover,
\[ ||u||_{p}\leq C_4,\]
(c) If $n\geq 3$, then for any $p\in (1, \frac{n(n+2)}{n^{2}-2})$ and $q\in (1, \frac{n+2}{n+1})$ 
the solution $u(x,t)$ has the following regularity in $Q_T$ for any $T>0$:
\[ u(x,t)\in L^{p}(Q_T)\bigcap L^q(0,T;W^{1,q}(\Omega)).\]
Moreover, there exists a constant $C_5$ such that
\[ ||u||_{L^{p}(Q_{T})}+||u||_{L^q(0,T:W^{1,q}(\Omega))}\leq C_5,\]
where $C_3, C_4$ and $C_5$ depend on $||f||_{L^{\infty}(0,\infty;L^{1}(\Omega))}, ||u_0||_{W^{1,q}(\Omega)}$, the upper bound of $T$ and known data. In addition,
$C_5$ also depends on $p$.

For $n=1$ and $n=2$, the conclusions of (a) and (b) are well-known from \cite{LSU}. For $n\geq 3$, the conclusion (c) is proved in \cite{BG1989}.

\hfill Q.E.D.

\ \\
{\bf Lemma 3.4.} Under the assumption H(2.1)-H(2.3) there exist two constants $r>0$ and $C_6$ such that
\[ \sup_{0<t<T}\int_{\Omega} u_1^{1+r} dx+ \int_{0}^{T}\int_{\Omega}|\nabla u_1^{\frac{1+r}{2}}|^2dxdt
\leq C_6,\]
where $r$ and $C_6$ depends only on known data.\\
{\bf Proof.} We assume that $n\geq 3$ since we already obtained a much stronger estimate for $n\leq 2$
by Lemma 3.3. We derive the estimate by using the standard energy method. Indeed, let $r>0$ to be chosen. We multiply
Eq.(2.1) by $u_1^r$ and integrate over $\Omega$ to obtain:
\begin{eqnarray*}
 & &  \frac{1}{1+r} \frac{d}{dt}\int_{\Omega}u_1^{1+r}dx+ \frac{4rd_0}{(1+r)^2}\int_{\Omega}|\nabla u_1^{\frac{1+r}{2}}|^2dx+d\int_{\Omega}u_1^{1+r}dx \\
& & \leq \sigma\int_{\Omega}u_3 u_1^r dx\\
& & \leq \varepsilon\int_{\Omega}u_1^{r\frac{q}{q-1}}dx+C(\varepsilon)\int_{\Omega}u_3^{q} dx.
\end{eqnarray*}
where at the final step we have used Young's inequality with a small parameter $\varepsilon >0$.

Now we choose $r=q-1>0$ and $\varepsilon =\frac{d}{2}$ to obtain
\begin{eqnarray*}
& & \frac{1}{1+r} \frac{d}{dt}\int_{\Omega}u_1^{1+r}dx+\frac{4rd_0}{(1+r)^2}\int_{\Omega}|\nabla u_1^{\frac{1+r}{2}}|^2dx+\frac{d}{2}\int_{\Omega}u_1^{1+r}dx \\
& & \leq C \int_{\Omega}u_3^{r+1} dx=C\int_{\Omega}u_{3}^{r+1}dx,
\end{eqnarray*}
which is uniformly bounded by Lemma 2.2 and 2.3 as long as $r\in (0,\frac{n(n+2)}{n^{2}-2}-1)$. Thus, the proof of Lemma 3.4 is completed.

\hfill Q.E.D.

With the estimates from previous lemmas, we are now ready to prove Theorem 2.1.

\ \\
{\bf Proof of Theorem 2.1}. 
There are several ways to prove the existence of a weak solution for the system (2.1)-(2.6).
 For any $\z>0$, we denote by $\chi(u)$ the Heaviside function. Set
 \[ \chi_{\z}(u)=\chi(u-\frac{1}{\z}).\]
\begin{eqnarray*}
f_{1\z}(Z) & = & b(x,t)-\beta_1u_{1}u_{2}-\beta_2 u_1h(B)-du_1+\sigma u_3(1-\chi_{\z}(u_3)),
\end{eqnarray*}
$f_2(Z), f_3(Z)$ and $f_4(Z)$ are the same as before.

Now we consider the following approximated reaction-diffusion system:
\setcounter{local}{1}
 \begin{eqnarray}
 & & u_{1 t}-\nabla[d_1\nabla u_{1}]=f_{1\z}(Z), \h (x,t)\in Q_{T},\\
 & & u_{2  t}-\nabla[d_2\nabla u_{2}]=f_{2}(Z), \h (x,t)\in Q_{T},\scl \\
& &  u_{3 t}-\nabla[d_3\nabla u_{3}]=f_3(Z), \h (x,t)\in Q_{T},\scl\\
 & & u_{4  t}-\nabla[d_4\nabla u_{4 }]=f_4(Z), \h (x,t)\in Q_T,\scl 
 \end{eqnarray}
subject to the same initial and boundary conditions (2.5)-(2.6).

We claim that the above approximate system has a unique weak solution in $V_2^4\bigcap L^{\infty}(Q_{T})$ for every small $\z>0$ and $||u_3||_{L^{\infty}(Q_{T})}\leq \frac{1}{\z}$.
Indeed, for every sufficiently small $\z>0$, if $u_3>\frac{1}{\z}$ in $Q_T$, then $u_1(x,t)$ satisfies
\[u_{1 t}-\nabla[d_1\nabla u_{1}]=b-\beta_1u_1u_2-\beta_2u_1h(u_4)-du_1,\]
subject to the same initial and boundary conditions.
The maximum principle implies that
\[ ||u_1||_{L^{\infty}(Q_{T}}\leq C,\]
where $C$ depends only on known data.

From Eq.(3.2), by using the maximum principle we see that $u_2$ will be uniformly bounded in $Q_T$ with the bound which has the same dependency as $u_1$. Thus, by Eq.(3.3), we apply the maximum principle again to see that $u_3$ would be bounded by a constant $C$ which depends only on known data. This is a contradiction, which implies that
$u_3\leq \frac{1}{\z}$ in $Q_{T}$ for a sufficiently small $\z$.

Consequently, the above system (3.1)-(3.4) associated with initial and boundary conditions (2.5)-(2.5)
has a unique weak solution 
\[ Z_{\z}(x,t)=(u_{1\z}, u_{2\z}, u_{3\z}, u_{4\z})\in V_{2}(Q_{T})^4\bigcap L^{\infty}(Q_{T}).\] 
(see \cite{EVANS} or \cite{LSU} ).
Moreover,
\[ 0\leq u_{3\z}\leq \frac{1}{\z},\h (x,t)\in Q_{T}.\]
Furthermore, it is clear that all apriori estimates from Lemma 3.1 to 3.4 hold and these bounds are independent  of $\z$. 
Note that $Z_{\z t}\in L^{p'}(0,T;W^{1,-p'}(\Omega))$ with $p'=\frac{1+r}{r}$.  By using the weak compactness of $L^{1+r}(Q_T)$, we can extract a subsequence of $\z$ if necessary,  as $\z\rightarrow 0$, that for $k=1,2,3,4,$
\begin{eqnarray*}
& & u_{k\z}(x,t)\rightarrow u_k(x,t), \h \mbox{strongly in $L^{1+r}(Q_{T})$},\\
& & \nabla u_{k\z} (x,t) \rightarrow \nabla u_k(x,t), \mbox{weakly in $L^{1+r}(Q_{T})$},\\
& & u_{k\z}(x,t)\rightarrow u_k(x,t), \h \mbox{a.e. in $Q_{T}$}.
\end{eqnarray*}
Moreover, by using Egorov's theorem, we see, as $\z\rightarrow 0$,
\[ \int\int_{Q_{T}}|u_{1\z}u_{2\z}-u_1u_2|dxdt\leq \int\int_{Q_{T}}[|(u_{1\z}-u_1)u_{2\z}|+u_1|u_{2\z}-u_2|]dxdt\rightarrow 0.\]

On the other hand, since $u_3\in L^{\infty}(0,T;W^{1,p}(\Omega))$, we see
\[ |\{(x,t)\in Q_{T}: u_3>\frac{1}{\z}\}|\]
is sufficiently small, provided that $\z$ is chosen to be sufficiently small.
It follows that for any test function $\phi(x,t)\in L^2(0,T;W^{1,p'}(\Omega)) $
\[\lim_{\z\rightarrow 0} \int\int_{Q_{T}}u_{3\z}(1-\chi_{\z}(u_{3\z})\phi dxdt= \int\int_{\Omega}u_{3}\phi dxdt.\]
Consequently, for all $k$ we have, as $\z \rightarrow 0$,
\[ f_{k\z}(Z_{\z})\rightarrow f_{k}(Z)\, a.e.\mbox{ in $Q_{T}$ and strongly in $L^1(Q_{T})$}.\]
For any test function $\phi_k\in X^*$ with $\phi_k\in L^{\infty}(Q_{T})$ and $\phi_{kt}\in X^*, \phi_k(x,T)=0$, we have for all $k=1,2,3,4$,
\begin{eqnarray*}
 & & \int_{0}^{T}\int_{\Omega}[-u_{k\z}\cdot \phi_{kt}+d_k\nabla u_{\z}\cdot \nabla\phi_k]dxdt=\int_{\Omega}u_{k0}(x)\phi_k(x,0)dx+\int_{\Omega} f_{\z k}(Z_{\z})\phi _k(x,t)dxdt. 
 \end{eqnarray*}
 After taking limit  as $\z\rightarrow 0$, we obtain a weak solution $Z(x,t)\in X$.
Moreover, by Lemma 3.4, we see $u_3\in W_{p}^{2,1}(Q_{T}).$

\hfill Q.E.D.

\ \\
{\bf Proof of Corollary 2.1}. Leq $n\leq 2.$
To see the uniform boundedness of $Z(x,t)$ in $Q$, we note that
$u_3(x,t)\in L^{\infty}(0, \infty; L^{1+r}(\Omega))$ for any $r\in (0, \frac{2n}{n+1})$. We apply a result of Lemma 2.6 from \cite{CC2017} that
\[ \sup_{0<t<\infty}||u_3||_{L^{\infty}(\Omega)} \leq C.\]

Now from Eq.(2.1) we see
\[ u_{1t}-\nabla [d_1 \nabla u_1]+du_1\leq b+\sigma u_3\]
subject to a homogeneous Neumann boundary condition and an initial condition.
Hence, since $u_3$ is uniformly bounded in $Q$ we see by using the maximum principle that
\[ \sup_{0<t<\infty}||u_1||_{L^{\infty}(\Omega)} \leq C,\]
where $C$ depends only on known data.

To prove the uniform boundedness for $u_2$, we use the same energy method as for $u_1$ to obtain
\begin{eqnarray*}
& & \frac{1}{1+r} \frac{d}{dt}\int_{\Omega}u_2^{1+r}dx+ \frac{4r d_0}{(1+r)^2}\int_{\Omega}|\nabla u_2^{\frac{1+r}{2}}|^2dx+(d+\gamma)\int_{\Omega}u_2^{1+r}dx \\
& & \leq C\int_{\Omega}u_2^{1+r}dx+C
\end{eqnarray*}
where we have used the uniform boundedness of $u_1$ and $C$ depends only on known data.

We recall a Gagliardo-Nirenberg's inequality:
 \[ ||u||_{L^{p}(\Omega)}\leq C||\nabla u||_{L^{q}(\Omega)}^{\theta}||u||_{L^{s}(\Omega)}^{1-\theta}+C||u||_{L^{1}(\Omega)},\]
 where $p,q, s$ and $\theta$ satisfy
 \[ \frac{1}{p}=\theta(\frac{1}{q}-\frac{1}{n})+(1-\theta)\frac{1}{s}, 0\leq \theta<1, 1\leq p<\infty.\]
 
For $n=2$, we choose $q=2, s=1, \theta=\frac{1}{2} $ and $p=2$,
 \[ ||u_2||_{L^{2}(\Omega)}\leq C||\nabla u_2||_{L^{2}(\Omega)}^{\theta}+C,\]
 since $||u_2||_{L^{1}(\Omega)}$ is uniformly bounded.
 
 That is
 \begin{eqnarray*}
 \int_{\Omega}u_{2}^{2}dx & \leq & C\left(\int_{\Omega}|\nabla u_2|^2dx\right)^{\frac{1}{2}}+C\\
 & \leq & \varepsilon \int_{\Omega}|\nabla u_2|^2dx +C(\varepsilon),
 \end{eqnarray*}
 
 If we choose $\varepsilon=\frac{rd_{0}}{(1+r)^2}$, we immediately obtain
 \[ y_2^{'}(t)+y_2(t)\leq C,\]
 where
 \[ y_2(t)=\frac{1}{2}\int_{\Omega}u_{2}(x,t)^{2}dx,\]
 $C$ depends only on known data. 

Hence,
\[ \sup_{0<t<\infty} y_2(t)\leq C, \]
where $C$ depends only on known data.

As $u_1$ and $u_2$ are uniformly bounded, we see
\[ \sup_{0<t<\infty}||f_2(z)||_{L^{2}(\Omega)} \leq C,\]
where $C$ depends only on known data.

Consequently, for $n=2$ we use Lemma 2.6 in \cite{CC2017} to have
 \[ ||u_2||_{L^{\infty}(Q)}\leq C,\]
 where $C$ depends only on known data.
 
 Once $u_2$ is uniformly bounded in $Q$, from Eq.(2.4) we can apply the maximum principle to obtain 
 \[  ||u_4||_{L^{\infty}(Q)}\leq C,\]
 where $C$ depends only on known data.

Moreover, when $f_k(Z)$ is bounded for $k=1,2,3,4,$, then the H\"older continuity of the weak solution $Z(x,t)$ directly comes from the standard DiGorgi-Nash's estimate for parabolic equations (\cite{LSU}).
The uniqueness of weak solution for $n\leq 2$ is straight forward by the energy method since $u_k$ for all $k$ is bounded in $Q$. 
 Hence, the proof of Corollary 2.1 is completed.
 
\hfill Q.E.D.

To prove Theorem 2.2, we need additional regularity conditions for $d_3(x,t)$ and other known data.
In the rest of this section we assume H(2.4) holds and $n = 3$.

Introduce
\[ Z_i(x,t):=\frac{\partial}{\partial x_{i}}(u_1(x,t),u_2(x,t),u_3(x,t),u_4(x,t)),i=1,2,\cdots n.\]

\ \\
{\bf Lemma 3.5.} Let the assumptions H(2.1)-H(2.4) hold. For any $q\in (1, \frac{n+2}{n+1})$,
there exists a constant $C_6$ such that
\[ \sum_{k=1}^{n}|| u_{3x_{k}}||_{W_{q}^{2,1}(Q_{T})}\leq C_6,\]
where $C_6$ depends only on known data.\\
{\bf Proof.} First of all, we note that for any $p\in (1, \frac{n+2}{n+1})$
\[||f_3||_{L^{p}(Q_{T})}\leq C,\]
where $C$ depends on known data and upper bound of $T$.

By $W_{p}^{2,1}(Q_{T})$-estimate, we obtain
\[ ||u_3||_{W_{p}^{2,1}(Q_{T})}\leq C.\]

Set $U_{ik}=u_{ix_{k}}, i=1,2,3,4, k=1, \cdots, n.$

We first take derivative for Eq.(2.3) with respect to $x_k$ to see that $U_{3k}$ satisfies
\[ U_{3kt}-\nabla [d_3 \nabla  U_{3k}] =\nabla [ d_{3x_{k}}\nabla u_{3}]+\gamma U_{2k}-U_{3k}(d+\sigma),\]

Since $U_{2k}\in L^{q}(Q_T)$ and $u_{3}\in W_{p}^{2,1}(Q_{T})$ by Lemma 3.4, we see
\[ ||\Delta u_3||_{p}+||U_{2k}||_{p}\leq C,\]
where $C$ depends only on known data.

Hence, we immediately obtain the interior $W_{p}^{2,1}(Q_T)$-estimate for $U_{3k}$ (\cite{Y1997}). Namely, for any 
$\Omega_0\subset  \Omega$ and $Q_{0T}=\Omega_0\times (0,T]$, there exists a constant $C$ such that
\[ ||U_{3k}||_{W_{p}^{2,1}(Q_{0T})}\leq C,\]
where $C$ depends only on the known constants and the distance $\Omega_0$ and $\partial \Omega$.

To derive the global $W_{p}^{2,1}$-estimate for $U_{3k}$, we first assume that 
$x_0\in \partial \Omega$ and $\partial \Omega$ is flat near $x_0$. Let
$B_r(x_0)$ be a small ball centered at $x_0$ with radius $r$ and 
\[ Q_{rT}(x_0)=(\Omega\bigcap B_{r}(x_0))\times (0,T], \Gamma_0=\{\partial \Omega \bigcap B_r(x_0): x_n=0.\} \]
The boundary condition for $u_{3}$ on $\Gamma_0$ becomes
\[ u_{3x_{n}}=0, \h (x,t)\in \Gamma_0\times (0,T].\]
It follows that, for $k=1, 2, \cdots, n-1,$ 
\[ u_{3x_{k}x_{n}}(x,t)=0, \h (x,t)\in \Gamma_0\times (0,T].\]
That is, for $k=1, 2, \cdots, n-1,$
\[  U_{3kx_{n}}(x,t)=0, \h (x,t)\in \Gamma_0\times (0,T].\]
Moreover, for $k=n$ we see
\[U_{3k}(x,t)=0, \h (x,t)\in \Gamma_0\times (0,T].\]
Now we can apply the global $W_{p}^{2,1}$-estimate for $U_{3k}$ for $k=1, 2, \cdots, n$ to obtain 
\[ ||U_{3k}||_{W_{q}^{2,1}(B_{r}(Q_{rT}(x_{0}))}\leq C,\]
where $C$ depends only on known data, $r$ and the upper bound of $T$.

When $\partial \Omega$ is not flat near $x_0$, since $\partial \Omega \in C^{2+\alpha}$ we use the standard transformation to convert $\Gamma_0$ to be a flat boundary near $x_0$ in a new coordinate and then use the same argument as above to obtain the desired $W_{q}^{2,1}$-estimate. We shall skip the step. The reader can find the detailed calculation for a general elliptic equation in \cite{EVANS}.
Finally, since $\partial \Omega$ is compact, after using a finite number of covering we can obtain the $W_{q}^{2,1}$-estimate near $\partial \Omega$. 

Thus, the proof for Lemma 3.5 is completed.

\hfill Q.E.D.

By using the standard embedding theorem for Sobolev spaces (see \cite{EVANS}), we immediately obtain the following consequence.
 \ \\
 {\bf Corollary 3.6. } Let $n\leq 3$. There exists a constant $C_7$ such that
 \[ ||u_3||_{C^{\alpha, \frac{\alpha}{2}}(\bar{Q}_T)}\leq C_7,\]
 where $q$ is the number in Lemma 3.5, $\alpha=\frac{n(q-1)}{n+2} \in (0,1)$
 and $C_7$ depends only on known data.

\ \\
{\bf Lemma 3.7.} Let $n\leq 3$. There exists a constant $C_8$ such that
\[ ||u_1||_{L^{\infty}(Q_{T})}\leq C_8.\]
where $C_8$ depends only on known data.\\
{\bf Proof}: Note that $u_1(x,t), u_2(x,t)\geq 0$. It follows that
\[ u_{1t}-\nabla[d_1 \nabla u_1]+du_1\leq b+u_3.\]
Since $u_3$ is bounded in $Q_T$ by Corollary 3.6, we apply the standard the comparison principle (\cite{LSU}) to obtain
\[ ||u_1||_{L^{\infty}(Q_{T})}\leq C_8,\]
where $C_8$ depends only on known data and the upper bound of $T$.

\ \\
{\bf Lemma 3.8.} Let $n\leq 3$. There exists a constant $C_9$ such that
\[ ||u_2||_{L^{\infty}(Q_{T})}+||u_4||_{L^{\infty}(Q_{T})}\leq C_9,\]
where $C_9$ depends only on known data.\\
{\bf Proof}: Since $u_1$ is uniformly bounded in $Q_T$ and $0\leq h(B)\leq 1$, we apply the standard the maximum principle to conclude that $ u_2$ is uniformly bounded:
\[ ||u_2||_{L^{\infty}(Q_{T})}\leq C[||I_0||_{L^{\infty}(\Omega)}+||u_1||_{L^{\infty}(Q_{T})}],\]
which is uniformly bounded by Lemma 3.7.

 When $u_2$ is uniformly bounded, we can use the maximum principle method to obtain an $L^{\infty}(Q_{T})$-bound for $u_4$.
 
 \hfill Q.E.D.

With the $L^{\infty}$-bounds for $u_k, k=1,2,3,4$ in hand, we can use DiGiorgi-Nash's estimate to obtain the following estimate in H\"older space (\cite{LSU}).

\ \\
{\bf Lemma 3.9.} Under the assumptions H(2.1)-H(2.4) if $n\leq 3$ there exists a constant $C_{10}$ 
and a number $\alpha\in (0,1)$ such that
\[ \sum_{k=1}^{4}||u_k||_{C^{\alpha, \frac{\alpha}{2}}(\bar{Q}_{T})}\leq C_{10},\]
where $C_{10}$ depends only on known data.

\hfill Q.E.D.

\ \\
{\bf Proof of Theorem 2.2}: With the apriori estimate in Lemma 3.9, we can prove the global existence by using either a bootstrap method (\cite{YCW2017}) or a fixed-point method. We use  Leray-Schauder's fixed point theorem.  Choose a Banach space $Y=L^{\infty}(Q_{T})$. For any $(\lambda, J)\in [0,1]\times Y,$ we consider the following reaction-diffusion system:
\setcounter{section}{3}
\setcounter{local}{5}
 \begin{eqnarray}
 u_{1t}-\nabla[d_1\nabla u_1] & = & \lambda b(x,t)-\beta_1Ju_2-\beta_2u_1 \cdot h(u_4)-du_1+\sigma u_3,\scl\\
 u_{2t}-\nabla[d_2\nabla u_2] & = & \beta_1J u_2+\beta_2u_1 \cdot h(u_4)-(d+\gamma)u_2,\scl \\
 u_{3t}-\nabla[d_r\nabla u_3] & = &  \lambda \gamma J(x,t)-(d+\sigma)u_3,\scl \\
 u_{4t}-\nabla[d_b \nabla u_4] & = & \xi u_2 +gu_4(1-\frac{u_4}{u_4+K})-\delta u_5. \scl
 \end{eqnarray}
 subject to the following initial and boundary conditions:
 \begin{eqnarray}
 & & \nabla_{\nu}Z(x,t) = 0,  (x,t)\in Q_{T},\scl\\
 & & Z(x,0)=\lambda (S_0(x),I_0(x),R_0(x),B_0(x)), \h x\in \Omega.\scl
 \end{eqnarray}
 
 For every $J(x,t)\in Y$ and $\lambda\in (0,1]$, under the assumption H(2.1)-(2.4) the standard theory for parabolic equation (\cite{LSU}) implies that the system (3.5)-(3.11) has a unique solution $Z^*(x,t)=(u_1^*(x,t),u_2^*(x,t),u_3^*(x,t),u_4^*(x,t))\in V_2(Q_{T})^4$. Moreover, the solution is H\"older continuous in $\bar{Q}_T$.
\[ ||Z^*||_{C^{\alpha,\frac{\alpha}{2}}(\bar{Q}_{T})} \leq C\]
 For every fixed $\lambda\in [0,1]$, we define
 a mapping $M_{\lambda}$ from  $ Y $ into $ Y$:
 \[ M_{\lambda}[J]=u_2^*(x,t).\]
 Note that $M_0[J]=0$ and a fixed-point of $M_1$ along with $u_1^*,u_3^*,u_4^*$ forms the solution of the original system (2.1)-(2.6).
 It is a routine to show that $M_{\lambda}$ is a continuous mapping from $Y$ into $Y$. Moreover,
 since  the embedding operator from $C^{\alpha,\frac{\alpha}{2}}(\bar{Q}_{T})$ into $Y$ is compact, the estimate in Lemma 3.9 implied that $M_{\lambda}$ is a compact mapping from $Y$ into $Y$. Furthermore, it is clear that all estimates from Lemma 3.6-3.9 hold for all fixed-points of $M_{\lambda}[u_2]=u_2$ and any $\lambda\in [0,1]$.
 By applying Leray-Schauder's fixed-point theorem, we see that $M_{\lambda}$ has a fixed-point.
 Particularly, when $\lambda=1$ the fixed point along with $u_1,u_3$ and $u_4$ forms a solution of the problem (2.1)-(2.6). The uniqueness is obvious since the solution is uniformly bounded in $Q_{T}$.

\hfill Q.E.D.

\section{Global Bounds and Asymptotic Behavior of Solution}

\setcounter{section}{4}
\setcounter{local}{1}
In this section we study the global asymptotic behavior of the solution and prove Theorem 2.3 and Theorem 2.4. In order to see a clear physical meaning for all species, we use the original variables
$S,I,R, B$ instead of $u_1,u_2,u_3, u_4$ in this section.

\ \\
{\bf Proof of Theorem 2.3}: Keep in mind that we shall always derive an estimate
for $S,I$ and $R$ as a first step. The second step for $B$ will be easy once we have the estimate for $I$ by Eq.(1.4).

We multiply Eq.(1.1) by $S$ and integrate over $\Omega$ to obtain
\begin{eqnarray*}
& & \frac{1}{2}\frac{d}{dt}\int_{\Omega}S^2dx+d_0\int_{\Omega}|\nabla S|^2dx+\beta_1\int_{\Omega}S^2Idx+\beta_2\int_{\Omega}S^2h(B)dx +d\int_{\Omega}S^2dx\\
& & \leq \int_{\Omega}b S dx+\sigma\int_{\Omega}SRdx\\
& & \leq  \int_{\Omega}b Sdx +\sigma\int_{\Omega}[\z R^2+\frac{1}{4\z}S^2]dx.
\end{eqnarray*}
where at the final step, we have used Cauchy-Schwarz's inequality with a small parameter $\z>0$.

Similarly, we mutiply Eq.(1.2) by  $ I$ and integrate of $\Omega$ to obtain
\begin{eqnarray*}
& & \frac{1}{2}\frac{d}{dt}\int_{\Omega}I^2dx+d_0\int_{\Omega}|\nabla I|^2dx+(d+\gamma)\int_{\Omega}I^2dx\\
& & \leq \beta_1\int_{\Omega}SI^2 dx+\beta_2\int_{\Omega}IS h(B)dx \\
& & \leq \beta_1\int_{\Omega}[\z S^2I+\frac{1}{4\z}I^2]dx+\beta_{2}\int_{\Omega}[\frac{1}{4\z}I^2+\z S^2h(B)]dx,
\end{eqnarray*}
where we have used the fact $0\leq h(B)\leq 1$.

We perform the same calculation for Eq.(1.3) to obtain
\begin{eqnarray*}
& & \frac{1}{2}\frac{d}{dt}\int_{\Omega}R^2dx+d_0\int_{\Omega}|\nabla R|^2dx+(d+\sigma)\int_{\Omega}R^2dx \\
& & \leq \gamma \int_{\Omega}IRdx\\
& & \leq \gamma \int_{\Omega}[\z I^2+\frac{1}{4\z}R^2]dx.
\end{eqnarray*}

Now we choose $\z=1$ and add up the above estimates to obtain
\begin{eqnarray*}
& & \frac{1}{2}\frac{d}{dt}\int_{\Omega}[S^2+I^2+R^2]dx+d_0\int_{\Omega}[|\nabla S|^2+\nabla I|^2+|\nabla R|^2]dx+ d\int_{\Omega} [S^2+I^2+R^2]dx\\
& & \leq \int_{\Omega}b S dx +\frac{\sigma}{4}\int_{\Omega}S^2dx+
\frac{\beta_1+\beta_2}{4}\int_{\Omega}I^2dx+\frac{\gamma}{4}\int_{\Omega} R^2dx\\
& & \leq \int_{\Omega}b Sdx +
\frac{\sigma+\beta_1+\beta_2+\gamma}{4}\int_{\Omega}[S^2+I^2+ R^2]dx.
\end{eqnarray*}

Let 
\[ Y(t)=\frac{1}{2}\int_{\Omega}[S^2+I^2+R^2]dx,\h  t\ge 0.\]
We see
\begin{eqnarray*}
& & Y'(t)+2(d-g_0)Y(t)  \leq 2\int_{\Omega} b S dx\\
& & \leq (d-g_0)\int_{\Omega} S^2dx +\frac{1}{(d-g_0)}\int_{\Omega} b^2 dx
\end{eqnarray*}
It follows that
\[ Y(t)\leq e^{-\frac{(d-g_0)t}{2}}[Y(0)-1]+\frac{2b_0|\Omega|}{(d-g_0)}, t\in [0,\infty).\]

The above estimate is good enough to derive the $L^{\infty}$-bound in $Q_T$ for any fixed $T>0$ by using an iteration technique and Sobolev embedding. Indeed, we use the same technique as in Lemma 3.5 in section 3 to obtain that $R_{x_i}\in W_{p_0}^{2,1}(Q_{T})$ for all $i=1, \cdots, n$, with $p_0>2$ and
\[ \sum_{i=1}^{n}||R_{x_i}||_p\leq C,\]
where $C$ depends only on known data and the upper bound of $T$.

By the Sobolev embedding, we see $R(x,t)\in L^{q_0}(Q_{T}$ with $q_0=\frac{(n+2)p_0}{n+2-4p_0}$ and
\[ ||R||_{q_0,Q_{T}}\leq C.\]
Hence, the energy estimate yields that $S(x,t)\in L^{q_0}(Q_{T})$ and
\[ ||S||_{q_0,Q_{T}}\leq C.\]
Next, we take derivative with respect to $x_i$ for Eq.(1.1)-(1.2) and follow the same calculations as for $R$, we are able to obtain $I_{i}\in L^{p_0}(Q_{T})$ for all $i$ and
\[  \sum_{i=1}^{n}||I_{x_i}||_{p_0,Q_{T}}\leq C,\]
where $C$ depends only on known data and $T$.

 Since all coefficients $d_i$ and $b(x,t), Z(,0)$ are smooth, after a finite number  of $l$-steps, we
 conclude that $R\in W_{p}^{2l,l}(Q_{T})$ and
 \[ ||R||_{W_{p}^{2l,l}(Q_{T})}\leq C,\]
 where $C$ depends only on known data.
 
 Consequently, if $l\geq \frac{n+2}{2}$ Sobolev's embedding yields
 \[ ||R||_{L^{\infty}(Q_{T})}\leq C,\]
 where $C$ depends only on known data and $T$.
 
 With the above $L^{\infty}(Q_{T})$-bound for $R(x,t)$, we immediately obtain an $L^{\infty}(Q_{T})$ bound for $S(x,t)$ from Eq.(1.1) by applying the maximum principle. Then from Eq.(1.2), the same maximum principle yields an $L^{\infty}(Q_{T})$-bound for $I$. Finally, from Eq.(1.4) we obtain an $L^{\infty}(Q_T)$ of $B(x,t)$ by applying a comparison principle.

Now for we use Lemma 2.6 to obtain
\[ \sum_{i=1}^{3}\sup_{Q}|u_i(x,t)|\leq C,\]
provided that $n\leq 4$, where $C$ depends only on known data.

The uniqueness of the weak solution can be proved easily since the weak solution is bounded.  

\hfill Q.E.D.

With the global bound for $Z(x,t)$ we are ready to prove the asymptotic behavior of the solution for the system (1.1)-(1.6).

\ \\
{\bf Proof of Theorem 2.4}: First of all, for we use Lemma 2.6 to obtain
\[ \sum_{i=1}^{3}\sup_{Q}|u_i(x,t)|\leq C,\]
provided that $n\leq 4$, where $C$ depends only on known data.

Let $S_{\infty}(x)$ be the following steady-state solution of the elliptic equation:
\begin{eqnarray}
& & -\nabla[d_1(x)\nabla S_{\infty}]=b_0(x)-dS_{\infty}, \h x\in \Omega,\\
& & \nabla_{\nu}S_{\infty}=0, \h x\in \partial \Omega.\scl
\end{eqnarray}
Obviously, the elliptic problem (4.1)-(4.2) has a unique solution (\cite{EVANS}). Moreover, since $d$ is positive, we can apply the maximum principle to see that $S_{\infty}(x)$ is uniformly bounded in $\Omega$. Furthermore, by Campanato estimate (see \cite{Y1997}) we have
\[||S_{\infty}||_{C^{1+\alpha}(\bar{\Omega})}\leq C,\]
where $C$ depends only on known data.

Set
\[ S^*(x,t)=S(x,t)-S_{\infty}(x), (I^*(x,t), R^*(x,t), B^*(x,t))=(I(x,t),R(x,t),B(x,t)), (x,t)\in Q.\]
We also define
\[ d_{i}^{*}(x,t)=d_i(x,t)-d_{i\infty}(x), \h i=1, 2,3,4.\]
Define
\begin{eqnarray*}
& & J_1(t)=\frac{1}{2}\int_{\Omega} S^*(x,t)^2dx, J_2(t)=\frac{1}{2}\int_{\Omega} I^*(x,t)^2dx, \\
& & J_3(t)=\frac{1}{2}\int_{\Omega} R^*(x,t)^2dx,J_4(t)=\frac{1}{2}\int_{\Omega} B^*(x,t)^2dx.
\end{eqnarray*}

Now, it is clear that $S^*(x,t)$ satisfies 
\begin{eqnarray}
& & S_t^*-\nabla[d_1(x,t)\nabla S^*]-\nabla[(d_1(x,t)-d_{1\infty}(x))\nabla S_{\infty}]\nonumber\\
& & = [b(x,t)-b_{0}(x)]-\beta_1SI-\beta_2Sh(B)-dS^*+\sigma R,\h (x,t)\in Q,\scl\\
& & \nabla_{\nu}S^*(x,t)=0, \h (x,t)\in \partial \Omega \times (0,\infty),\scl\\
& & S^*(x,0)=S_0(x)-S_{\infty}(x), \h x\in \Omega.\scl
\end{eqnarray}
For $I^*,R^*$ and $B^*$, they satisfy the same equations as $I,R,B$.

We multiply Eq.(4.3) by $S^*$ and integrate over $\Omega$ to obtain
\begin{eqnarray*}
& & J_1'(t)+d_0\int_{\Omega}|\nabla S^*|^2dx +dJ_1(t)\\
& &  \leq m_0\int_{\Omega}[|d_{1}^{*}| \cdot |\nabla S^*|]dx+\int_{\Omega}|b^*| |S^*| dx+\sigma\int_{Q}|R S^* | dx,
\end{eqnarray*}
where  $m_0=\sup_{\Omega}|\nabla S_{\infty}(x)|$.

We use Cauchy-Schwarz's inequality with a small parameter $\z$ to obtain
\begin{eqnarray*}
\int_{\Omega}[|d_1^{*}| |\nabla S^*|] dx & \leq &\varepsilon \int_{\Omega} |\nabla S^*|^2 dx+\frac{1}{4\varepsilon} \int_{\Omega} d_{1}^{*2}dx;\\
\end{eqnarray*}
By choosing $\varepsilon =\frac{d_0}{2m_{0}}$, we see
\begin{eqnarray*}
J_1'(t)+ +d J_1(t)
& &  \leq C\int_{\Omega}(d_{1}^{*})^2dx+\int_{\Omega}[|b_1^*||S^*]dx+\sigma \int_{\Omega}R^2 dx,\\
& & \leq C(\z_0)\int_{\Omega}[|d_{1}^{*}|^2+|b_1^*|^2]dx +\z_0\int_{\Omega}|S^*|^2 dx+\sigma \int_{\Omega}R^2 dx
\end{eqnarray*}
where $\z_0>0$ is any small number and $C(\z_0)$ depends only on known data and $\z_0$.

On the other hand, note that $I^*(x,t), R^*(x,t)$ and $B^*(x,t)$ satisfy the same equations as $I,R$ and $B$. If we perform the same energy estimate, then we have the same estimate as $I,R$ and $B$ except additional terms involving $d_k^*$:
\[ d_{2}^*(\nabla I^*\cdot \nabla I),\, d_{3}^*(\nabla R^*\cdot \nabla R).\]
we follow exactly the same method as for $J_1(t)$  to obtain the estimates similar to the one in the proof of Theorem 2.3. For these additional terms, we use the same argument as for $S^*$ to obtain
\begin{eqnarray*}
& & \int_{\Omega} |d_{2}^*\nabla I^*\cdot \nabla I|dx\leq \z\int_{\Omega}|\nabla I^*|^2sdxx+C(\z)
\int_{\Omega} |d_2^*|^2dx;\\
& & \int_{\Omega} |d_{3}^*\nabla R^*\cdot \nabla R|dx\leq \z\int_{\Omega}|\nabla R^*|^2sdx+C(\z)
\int_{\Omega} |d_3^*|^2dx;
\end{eqnarray*}
where $\z>0$ is an arbitrary small constant.

We choose $\z=\frac{d_0}{2}$. Let 
\[ J(t)=J_1(t)+J_2(t)+J_3(t), t\geq 0.\]
We  conclude
\[ J'(t)+(d-g_0-\z_0)J(t)\leq C\int_{\Omega}[\sum_{i=1}^{3}|d_i|^2+|b^*|^2]dx.\]
where 
\[ g_0=\frac{\sigma+\beta_1+\beta_2+\gamma}{4}\]
and $C$ depends only on known data.

We choose 
\[ \z_0=\frac{d-g_0}{2}>0.\]

Then
\[ J(t) \leq e^{-\frac{(d-g_0)t}{2}}[J(0)+C]+CM(t)],\]
where 
\[ M(t)=\sup_{0<\tau<t}\int_{\Omega}[\sum_{k=1}^{3}|d_k^*(x,\tau)|^2+|b^*(x,\tau)|^2]dx.\]

By the assumption H(2.5), we  obtain
\[ \lim_{t\rightarrow \infty}J(t)=0.\]

Once we know $I(x,t)$ goes to $0$ in $L^2(\Omega)$-sense as $t \rightarrow \infty$, we can use the same argument for $J_4(t)$ (see \cite{Temam}) to obtain
\[ \lim_{t\rightarrow \infty}J_4(t)=0.\]
By the definition of a global attractor (see \cite{Temam}), we see that $Z_{\infty}(x)=(S_{\infty}(x),0,0,0)$ is a global attractor.
Thus, the proof of Theorem 2.4 is completed.

\hfill Q.E.D.

 \section{Conclusion}
 
 In this paper we studied a mathematical model for the Cholera epidemic without life-time immunity.
 The model equations are governed by a coupled reaction-diffusion system with different diffusion coefficients  for each species. We established the global well-posedness for the coupled reaction-diffusion system under a certain condition on known data. Moreover, the long-time behavior of the solution is obtained for any space dimension $n$. Particularly, we prove that there is a global attractor for the system under appropriate conditions on known data. These results justify the mathematical model and provide scientists a deeper understanding of the dynamics of interaction between
 bacteria and susceptible, infected and recovered human hosts.
The mathematical model with the help of real data analysis provides a scientific foundation for policy-makers to make better decisions for the general public in health and medical sciences. The main tools used in this paper come from some delicate theories for elliptic and parabolic equations. The method developed in this paper can be used to study other models such as the avian influenza for birds.
 
 \ \\
 {\bf Acknowledgements}
 
   This paper was motivated by some open questions raised by Professor K. Yamazaki in a seminar at Washington State University. The author would like to thank Professor Yamazaki and Professor X. Wang
   for some interesting discussions. Many thanks also go to Mr. Brian Yin, Esq., from
   Law Firm Clifford Chance US LLP, who helped to edit the original paper.

\end{document}